\documentclass[twoside]{article}
\include{graphicx}
\title{On the asymptotics of a sequence of lacunary binomial-type polynomials}
\author{R. B. Paris\\
\\
{\em School of Engineering, Computing and Applied Mathematics,}\\
{\em University of Abertay Dundee, Dundee DD1 1HG, UK}\\
E-Mail: r.paris@abertay.ac.uk}
\begin{document}
\def\f#1#2{\mbox{${\textstyle \frac{#1}{#2}}$}}
\def\dfrac#1#2{\displaystyle{\frac{#1}{#2}}}
\def\boldal{\mbox{\boldmath $\alpha$}}
\newcommand{\bee}{\begin{equation}}
\newcommand{\ee}{\end{equation}}
\newcommand{\lam}{\lambda}
\newcommand{\ka}{\kappa}
\newcommand{\al}{\alpha}
\newcommand{\om}{\omega}
\newcommand{\Om}{\Omega}
\newcommand{\fr}{\frac{1}{2}}
\newcommand{\fs}{\f{1}{2}}
\newcommand{\g}{\Gamma}
\newcommand{\br}{\biggr}
\newcommand{\bl}{\biggl}
\newcommand{\ra}{\rightarrow}
\renewcommand{\topfraction}{0.9}
\renewcommand{\bottomfraction}{0.9}
\renewcommand{\textfraction}{0.05}
\newcommand{\mcol}{\multicolumn}
\date{}
\maketitle
\pagestyle{myheadings}
\markboth{\hfill \sc R. B.\ Paris \hfill}
{\hfill {\it Lacunary binomial-type polynomials} \hfill}
\begin{abstract}
We examine the asymptotics of a sequence of lacunary binomial-type polynomials $\wp_n(z)$ as $n\ra\infty$ that have arisen in the problem of
the expected number of independent sets of vertices of finite simple graphs. We extend the recent analysis of Gawronski and Neuschel by employing the method of steepest descents applied to an integral representation.  
The case of complex $z$ with $|z|<1$ is also considered. Numerical results are presented to illustrate the accuracy of the resulting expansions.
\vspace{0.4cm}

\noindent {\bf Mathematics Subject Classification:} 30E15, 41A60 
\vspace{0.3cm}

\noindent {\bf Keywords:} Asymptotics, lacunary polynomials, saddle points
\end{abstract}
\vspace{0.3cm}

\begin{center}
{\bf 1. \  Introduction}
\end{center}
\setcounter{section}{1}
\setcounter{equation}{0}
\renewcommand{\theequation}{\arabic{section}.\arabic{equation}}
The sequence of lacunary binomial-type polynomials defined by
\bee\label{e11}
\wp_n(z)=\sum_{k=0}^n \left(\!\!\begin{array}{c}n\\k\end{array}\!\!\right)z^{k(k-1)/2}
\ee
has been encountered by Brown {\it et al.\/}  \cite{BDM} in connection with the problem of the expected number of independent sets of vertices of finite simple graphs. These authors obtained upper and lower bounds for the values of $\wp_n(z)$ for $0<z<1$ and, based on analytic results and numerical computations, made the following conjecture on its asymptotic behaviour as $n\ra\infty$
\bee\label{e12}
\wp_n(1/y)\sim\frac{1}{\sqrt{w(n)}}\,\exp \left(\frac{w(n)^2+2w(n)}{2\log\,y}\right)
\ee
when $y>1$. Here, $w(n)=W(n\surd y \log\,y)$, where $W(a)$ denotes the Lambert $W$-function
which, for $a>0$, is defined as the positive solution of the equation $te^t=a$ \cite[p.~111]{DLMF}. 

Recently, Gawronski and Neuschel \cite{GN} considered an integral representation for $\wp_n(z)$
and employed a path displacement argument in the complex plane combined with a non-standard version of the saddle-point method to obtain the asymptotic formula
when $y>1$ 
\bee\label{e13}
\wp_n(1/y)=\frac{1}{\sqrt{r(n)}}\,\exp \left(\frac{r(n)^2+2r(n)}{2\log\,y}\right)\left\{\Theta(y)+o(1)\right\}
\ee
as $n\ra\infty$, where 
\bee\label{e13a}
\Theta(y)=1+2\sum_{k=1}^\infty e^{-2\pi^2k^2/\log\,y} \cos\left(\frac{2\pi kr(n)}{\log\,y}\right)
\ee
and the quantity $r(n)$ is defined as the positive solution of the equation
\[t(e^t+\surd y)=n\surd y\,\log\,y.\]
It was established in \cite{GN} that the expression in front of the curly braces in (\ref{e13})
is asymptotically equivalent to the conjectured approximation in (\ref{e12}) in the limit $n\ra\infty$.

In this paper, we consider the same integral representation for $\wp_n(z)$ as derived in \cite{GN} and employ the method of steepest descents to determine a more accurate expansion for moderately large 
values of $n$. We also discuss the case of complex $z$ satisfying $|z|<1$. We present numerical
computations to demonstrate the accuracy of our expansion and the approximation in (\ref{e13}).

\vspace{0.6cm}

\begin{center}
{\bf 2. \ An integral representation and saddle-point structure when $z<1$}
\end{center}
\setcounter{section}{2}
\setcounter{equation}{0}
\renewcommand{\theequation}{\arabic{section}.\arabic{equation}}
For convenience in presentation we set throughout $z=1/x^2$, where $|x|>1$ with $|\arg\,x|\leq\fs\pi$.
Substituting the standard result (cf. \cite[Eq.~(7.7.3)]{DLMF})
\[x^{-k^2}=\frac{1}{2\sqrt{\pi \log\,x}}\int_{-\infty}^\infty \exp \left(-\frac{s^2}{4\log\,x}+iks\right)ds\qquad (|x|>1)\]
into (\ref{e11}), followed by an interchange in the order of summation and integration and use of the binomial theorem, we obtain the integral representation derived in \cite{GN}
\begin{eqnarray}
\wp_n(x^{-2})&=&\frac{1}{2\sqrt{\pi\log\,x}}\int_{-\infty}^\infty \exp\left(\frac{-s^2}{4\log\,x}\right)(1+xe^{is})^nds\nonumber\\
&=&\frac{1}{2\sqrt{\pi\log\,x}}\int_{-\infty}^\infty e^{-n\psi(s)}ds,\label{e21}
\end{eqnarray}
where
\bee\label{e21a}
\psi(s)=\frac{s^2}{4n\log\,x}-\log(1+xe^{is}).
\ee
We proceed in Section 3 with the determination of the asymptotics of the integral in (\ref{e21}) for $n\ra\infty$ and finite real $x$ by application of the method of steepest descents.
In Section 4 we extend these considerations to complex $x$. Before doing this we first examine the saddle-point structure of the phase function $\psi(s)$ when $x$ is real and $x>1$.

The saddle points $s_k$ of $\psi(s)$ are given by $\psi'(s)=0$, or
\bee\label{e22}
s_ke^{-is_k}(1+xe^{is_k})=2inx\log\,x.
\ee
There is a single saddle situated on the imaginary axis together with two infinite strings of saddles symmetrically
positioned about $\mbox{Re} (s)=0$ resulting from the periodicity of the function $e^{is}$. If we write
\bee\label{e2v}
s_k=\sigma_k+2\pi k
\ee
for integer $k$ we have from (\ref{e22}), with $\alpha:=2x\log\,x$,
\[\log\,(\sigma_k+2\pi k)-i\sigma_k+\log\,(1+xe^{i\sigma_k})=\log\,(in\alpha).\]

A straightforward perturbative solution of this equation for large $n$ can be obtained by putting $\sigma_k=i\log\,n-\epsilon$, where $|\epsilon|\ll\log\,n$, to find
\[\epsilon\simeq i(\log\log\,n-\log\,\alpha)+i\log \left(1-\frac{i(2\pi k-\log\log\,n)}{\log\,n}\right).
\]
This yields the approximation valid for large $n$ and bounded $k$ given by
\bee\label{e23}
\sigma_k\simeq i\log\left(\frac{n\alpha}{\log\,n}\right)-i\log\,\left(1-\frac{i(2\pi k-\log\log\,n)}{\log\,n}\right).
\ee
The saddle on the imaginary axis corresponding to $k=0$ is 
\bee\label{e24}
s_0=\sigma_0=ir(n),\qquad r(n)\simeq \log\left(\frac{n\alpha}{\log\,n}\right)+ \frac{\log\log\,n}{\log\,n}.
\ee
Provided $2\pi|k|\ll \log\,n$, the saddles with $|k|\geq 1$ have from (\ref{e23}) the approximate location
\bee\label{e24a}
\sigma_k\simeq -\frac{2\pi k}{\log\,n}+i\left\{r(n)-\frac{1}{2}\left(\frac{2\pi k-\log\log\,n}{\log\,n}\right)^{\!\!2}\right\}.
\ee
This last result indicates that the first few saddles situated in the right and left half-planes have $\mbox{Im} (s_k)\simeq r(n)$
for $n\ra\infty$. However, the value of $n$ required for the validity of this approximation is extremely large.
When $k=\pm 1$, the requirement $2\pi |k|/\log\,n=10^{-1}$, for example, yields $n=e^{20\pi}\simeq 2\times 10^{27}$. For moderately large $n$ the values of $\mbox{Im} (s_k)$ for $|k|\geq 1$ are found to be significantly less than $r(n)$; see Table 1. 
\begin{table}[th]
\caption{\footnotesize{The location of the saddles $s_k$ for $0\leq k\leq 5$ and their approximate values from (\ref{e23}) when $n=1000$, $x=2$. The saddles $s_{-k}=-{\overline s}_k$, where the bar denotes the complex conjugate. }}
\begin{center}
\begin{tabular}{c|r|r}
\mcol{1}{c|}{$k$} & \mcol{1}{c|}{$s_k$} & \mcol{1}{c}{Approximate $s_k$}\\
[.05cm]\hline
&& \\[-0.2cm]
0  & $6.112742i$            &          $6.323089i$\\
1  & $5.521734+5.839316i$  & $5.382118+5.846300i$\\
2  & $11.427821+5.387286i$ & $11.372547+5.323709i$\\
3  & $17.544733+5.019893i$ & $17.536813+4.957375i$\\
4  & $23.741718+4.737505i$ & $23.757372+4.684147i$\\
5  & $29.972889+4.513009i$ & $30.002189+4.467837i$\\
[.15cm]\hline
\end{tabular}
\end{center}
\end{table}

The paths of steepest descent, which we denote by ${\cal C}_k$, and ascent are given by the paths on which
\[\mbox{Im}(\psi(s)-\psi(s_k))=0.\]
The steepest descent paths can terminate either at infinity in the sectors $|\arg\,s|<\f{1}{4}\pi$ or at the points where $(1+xe^{is})=0$; that is, at the logarithmic singularities of $\psi(s)$ given by
\bee\label{e25}
T_k=i \log\,x+(2k+1)\pi \qquad (k\in Z).
\ee
It is found (we omit these details) that as $k\ra\pm\infty$ the saddles $s_k$ approach the points $T_k$ in (\ref{e25}). 
Paths of steepest ascent terminate at infinity in the sectors $|\arg\,(\pm is)|<\f{1}{4}\pi$. An example of the distribution of the saddles and the paths of steepest descent and ascent is shown in Fig.~1. The case of complex $x$ is deferred to Section 4.
\begin{figure}[ht]
	\begin{center}
		\includegraphics[width=0.60\textwidth]{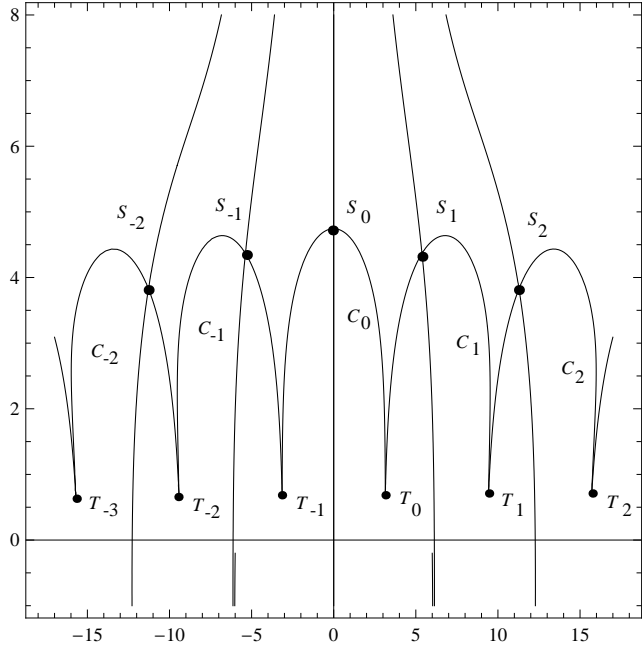}
	\caption{\small{The steepest descent paths ${\cal C}_k$ through the saddles $s_k$ (heavy dots) with $|k|\leq 2$ and the associated logarithmic singularities $T_k$ given by (\ref{e25}) when $n=200$, $x=2$. The steepest ascent paths pass to infinity parallel to the imaginary $s$-axis.}}\label{f1}
	\end{center}
\end{figure}
\vspace{0.6cm}

\begin{center}
{\bf 3. \ Derivation of the expansion for $x>1$}
\end{center}
\setcounter{section}{3}
\setcounter{equation}{0}
\renewcommand{\theequation}{\arabic{section}.\arabic{equation}}
By a standard application of Cauchy's theorem we can shift the integration path in (\ref{e21}) from the real $s$-axis to a parallel line through the points $T_k$ (where the integrand vanishes), and thence to coincide with the steepest descent paths ${\cal C}_k$ through the saddles $s_k$. Our final integration path then takes on a ``serpentine''
form as it passes over each saddle in turn. Hence we can write (when $x$ is real)
\bee\label{e29b}
\wp_n(x^{-2})=J_0+2\mbox{Re} \sum_{k=1}^\infty J_k,
\ee
where
\bee\label{e29a}
J_k=\frac{1}{2\sqrt{\pi \log\,x}}\int_{{\cal C}_k} e^{-n\psi(s)}ds.
\ee

We present only the calculation of the expansion of $J_0$, since the details for $J_k$ ($k\geq 1$) are similar. We define the quantities
\bee\label{e29}
\lambda_k:=xe^{is_k},\qquad a_k:=\frac{(1+\lambda_k)^2}{4n\lambda_k \log\,x},\qquad
\om_k:=\frac{2n\lambda_k\log\,x}{(1+\lambda_k)^2}\quad (k\geq 0),
\ee
so that 
\[\psi''(s_0)=\frac{1}{2n\log\,x}+\frac{\lambda_0}{(1+\lambda_0)^2}=\frac{1+\om_0}{2n\log\,x}.\]
Then application of the method of steepest descents \cite[p.~48]{DLMF} produces
\[J_0\sim \frac{e^{-n\psi(s_0)}}{\sqrt{2\pi\psi''(s_0)\log\,x}} \sum_{j=0}^\infty \frac{\g(j+\fs)\, d_{j0}}{n^{j+\fr}}\]
as $n\ra\infty$. The coefficients $d_{j0}$ ($0\leq j\leq 3$) are given by \cite[p.~119]{D}, \cite[p.~13]{P}
\[d_{00}=1,\qquad d_{10}=\frac{1}{12\psi''(s_0)}(5p_3^{\,2}-3p_4),\]
\[d_{20}=\frac{1}{864(\psi''(s_0))^2}(385p_3^{\,4}-630p_3^{\,2}p_4+168p_3p_5+105p_4^{\,2}-24p_6),\]
\[d_{30}=\frac{1}{777600(\psi''(s_0))^3}(425425p_3^{\,6}-1126125p_3^{\,4}p_4+675675p_3^{\,2}p_4^{\,2}-51975p_4^{\,3}+360360p_3^{\,3}p_5-
\]
\[\hspace{1cm}-249480p_3p_4p_5+13608p_5^{\,2}-83160p_3^{\,2}p_6+22680p_4p_6+12960p_3p_7-1080p_8), \]
where, for convenience in presentation, we have defined
\[p_r:=\frac{\psi^{(r)}(s_0)}{\psi''(s_0)}\qquad (r\geq 3).\]
Insertion of the values of the derivatives of $\psi(s)$ evaluated at $s_0$ then yields after some algebra with the help of {\it Mathematica} the coefficients in the form
\[d_{00}=1,\qquad d_{10}=\frac{-Q_1(a_0,\lambda_0)}{6(1+2a_0)^3 \lambda_0},\]
\[
d_{20}=\frac{Q_2(a_0,\lambda_0)}{216(1+2a_0)^6 \lambda_0^2},
\qquad
d_{30}=\frac{Q_3(a_0,\lambda_0)}{
97200(1+2a_0)^9 \lambda_0^3}. \]
The quantities $Q_j(a_0,\lambda_0)$ are defined in the Appendix in terms of polynomials in $\lambda_0$.

From (\ref{e22}), (\ref{e24}) and (\ref{e29}) it is easy to see that
\[\sigma_0\sim\log\,n,\quad e^{\sigma_0}\sim\frac{2nx\log\,x}{\log\,n},\quad \lambda_0\sim\frac{\log\,n}{2n\log\,x}\quad a_0\sim \frac{1}{2\log\,n}\qquad (n\ra\infty).\]
If we define the coefficients $c_{j0}$ by
\[c_{00}=1,\qquad c_{10}=\frac{-Q_1(a_0,\lambda_0)\chi}{6(1+2a_0)^3},\]
\bee\label{e2coeff}
c_{20}=\frac{Q_2(a_0,\lambda_0)\chi^2}{216(1+2a_0)^6},
\qquad
c_{30}=\frac{Q_3(a_0,\lambda_0)\chi^3}{
97200(1+2a_0)^9}, 
\ee
where we have put $\chi\equiv \log\,n/(n\lambda_0)$,
it follows from (\ref{a1}) that
\[c_{10}\sim -\f{1}{3} \log\,x,\qquad c_{20}\sim \f{1}{54} (\log\,x)^2,\qquad c_{30}\sim \f{139}{12150}(\log\,x)^3\]
as $n\ra\infty$, and therefore that $c_{j0}=O(1)$ ($1\leq j\leq 3$) in this limit.
Then we obtain finally
\bee\label{e30a}
J_0\sim\frac{e^{-n\psi(s_0)}}{\sqrt{1+\omega_0}}\sum_{j=0}\frac{(\fs)_j c_{j0}}{(\log\,n)^j}\qquad (n\ra\infty),
\ee
where $(\fs)_j=\g(j+\fs)/\g(\fs)$ is Pochhammer's symbol and $\sum_{j=0}$ means that the sum is restricted to the first few terms of the series.
The quantity $\omega_0 \sim \log\,n$ as $n\ra\infty$.
We remark that a full treatment would require a knowledge of the behaviour of the coefficients in general as $n\ra\infty$ to establish the asymptotic nature of (\ref{e30a}); we do not carry this out here.

The calculation of the integrals $J_k$ ($k\geq 1$) follows the same procedure and we find
\bee\label{e30b}
J_k\sim\frac{e^{-n\psi(s_k)}}{\sqrt{1+\omega_k}}\sum_{j=0}\frac{(\fs)_j c_{jk}}{(\log\,n)^j}\qquad (n\ra\infty),
\ee
where the first few coefficients $c_{jk}$ are defined in (\ref{e2coeff}) with $a_0$ and $\lambda_0$ replaced by $a_k$ and $\lambda_k$.
In terms of the quantities $\sigma_k$ defined in (\ref{e2v}), we have
\[\psi(s_k)=\psi(\sigma_k)+\frac{\pi^2k^2+\pi k\sigma_k}{n\log\,x}.\]
Then the expansion of $\wp_n(x^{-2})$ when $x>1$ is given by the following theorem.
\newtheorem{theorem}{Theorem}
\begin{theorem}
For real $x$ satisfying $x>1$ we have from (\ref{e29b}) the expansion
\[\wp_n(x^{-2})\sim\frac{e^{-n\psi(\sigma_0)}}{\sqrt{1+\omega_0}}\sum_{j=0}\frac{(\fs)_j c_{j0}}{(\log\,n)^j}\hspace{6cm}\]
\bee\label{e28}
\hspace{3cm}+2\mbox{Re} \sum_{k=1}^\infty \exp \left(-\frac{\pi^2k^2+\pi k\sigma_k}{\log\,x}\right)\,
\frac{e^{-n\psi(\sigma_k)}}{\sqrt{1+\omega_k}}\sum_{j=0}\frac{(\fs)_j c_{jk}}{(\log\,n)^j}
\ee
as $n\ra\infty$. The coefficients $c_{jk}$ and the quantities $\sigma_k$ and $\omega_k$ are defined in (\ref{e2coeff}), (\ref{e2v}) and (\ref{e29}), respectively.
\end{theorem}

From (\ref{e24a}), the saddles with finite $k$ approximately satisfy $\sigma_k\simeq ir(n)\simeq i\log\,n$
and $\om_k\simeq\log\,n\simeq r(n)$ as $n\ra\infty$. Since, with this approximation,
\[-\psi(\sigma_k)\simeq \frac{r(n)^2}{4n\log\,x}+\log \left(1+\frac{r(n)}{2n\log\,x}\right)\simeq\frac{r(n)^2+2r(n)}{4n\log\,x} \qquad (n\ra\infty),\]
it is seen that the leading terms of (\ref{e28}) agree with the result (\ref{e13}) (with $y=x^2$) obtained by Gawronski and Neuschel \cite{GN}.

In Table 2 we show the absolute relative error in the computation of $\wp_n(x^{-2})$ for different $n$ and $x$
using the asymptotic expansion (\ref{e28}) as a function of the truncation index $j$ in the dominant sum (corresponding to $k=0$). 
Table 3 compares the value of $\wp_n(x^{-2})$ with that obtained from (\ref{e28}) (with truncation index $j=3$ in the dominant sum) and the approximation (\ref{e13}). In each table
two terms have been employed in the sum with $k=1$, with the contributions from the sums with $k\geq 2$ being negligible on account of the $\exp (-\pi^2k^2/\log\,x)$ dependence. It is seen that for the values of $n$ chosen ($n\leq 1000$) the formula (\ref{e13}) gives only a gross approximation, whereas (\ref{e28}) yields very good agreement.
\begin{table}[th]
\caption{\footnotesize{The absolute relative error in the computation of $\wp_n(x^{-2})$ for different $n$ and $x$ 
using the asymptotic expansion (\ref{e28}) as a function of the truncation index $j$ in the dominant sum.}}
\begin{center}
\begin{tabular}{c|c|c|c}
\mcol{1}{c|}{$j$} & \mcol{1}{c|}{$n=200,\ x=1.10$} & \mcol{1}{c|}{$n=400,\ x=1.05$} & \mcol{1}{c}{$n=400,\ x=1.50$}\\
[.05cm]\hline
&&& \\[-0.3cm]
0  & $1.338\times 10^{-3}$ & $6.729\times 10^{-4}$ & $5.802\times 10^{-3}$\\
1  & $1.913\times 10^{-5}$ & $4.989\times 10^{-6}$ & $1.716\times 10^{-4}$\\
2  & $3.197\times 10^{-7}$ & $4.340\times 10^{-8}$ & $2.179\times 10^{-6}$\\
3  & $3.053\times 10^{-9}$ & $\ 2.287\times 10^{-10}$& $4.686\times 10^{-7}$\\
[.15cm]\hline
\end{tabular}
\end{center}
\end{table}
\begin{table}[th]
\caption{\footnotesize{The values of $\wp_n(x^{-2})$ compared with the asymptotic expansion (\ref{e28})
(with $j=3$ in the dominant sum and $j=2$ in the sum with $k=1$) and the approximation (\ref{e13}) for different values of $n$ and $x$. }}
\begin{center}
\begin{tabular}{c|c|c}
\mcol{1}{c|}{} & \mcol{1}{c|}{$n=200,\ x=1.20$} & \mcol{1}{c}{$n=200,\ x=2$}\\
[.05cm]\hline
&& \\[-0.3cm]
$\wp_n(x^{-2})$ & $8.562122063\times 10^9$ & $4.398555252\times 10^4$\\
Asymptotic & $8.562122013\times 10^9$ & $4.398536817\times 10^4$\\
Eq. (1.3)  & $7.864432769\times 10^9$ & $4.712945605\times 10^4$\\
[.15cm]\hline
\mcol{1}{c|}{} & \mcol{1}{c|}{$n=400,\ x=1.10$} & \mcol{1}{c}{$n=1000,\ x=2$}\\
[.05cm]\hline
&& \\[-0.3cm]
$\wp_n(x^{-2})$ & $9.488964463\times 10^{18}$ & $2.202064917\times 10^7$\\
Asymptotic & $9.488964461\times 10^{18}$ & $2.202060088\times 10^7$\\
Eq. (1.3)  & $7.418083490\times 10^{18}$ & $2.370080869\times 10^7$\\
[.15cm]\hline
\end{tabular}
\end{center}
\end{table}
\vspace{0.6cm}

\begin{center}
{\bf 4. \  The case of complex $x$}
\end{center}
\setcounter{section}{4}
\setcounter{equation}{0}
\renewcommand{\theequation}{\arabic{section}.\arabic{equation}}
We now consider $x$ to be a complex variable and write $x=|x| e^{i\theta}$, where it will be sufficient to restrict our attention to values of the phase $\theta$ satisfying $0\leq\theta\leq\fs\pi$. In terms of the original variable $z$ this corresponds to $0\leq\arg\,z\leq\pi$.
When $\theta>0$, it can be seen from (\ref{e23}) and (\ref{e25}) that both the saddles $s_k$ and the logarithmic singularities $T_k$ are displaced to the left\footnote{When $\theta<0$ the situation is reversed: the saddles and the points $T_k$ are displaced to the right.}; see Fig.~2. A more significant difference between the $\theta=0$
and $\theta>0$ cases, however, is the connectivity of the saddles.

When $\theta=0$, it is found that Im$(\psi(s_k))$ increases monotonically with $k$ and all the paths of steepest descent through the saddles $s_k$ have their endpoints at their associated logarithmic singularities $T_{k-1}$ and $T_k$; see Fig.~1. When $\theta>0$, the quantities
Im$(\psi(s_k))$ attain a maximum for a certain $k\geq 1$ and thereafter decrease monotonically with $k$. This results in a Stokes phenomenon occurring at an infinite sequence of $\theta$-values when a pair of adjacent saddles connects; see Table 4. 
\begin{table}[h]
\caption{\footnotesize{Values of the phase $\theta$ when the saddles $s_k$, $s_{k+1}$ (with $1\leq k\leq 5$) connect to produce a Stokes phenomenon.}}
\begin{center}
\begin{tabular}{c|c|c}
\mcol{1}{c|}{Saddles} & \mcol{1}{c|}{$n=200,\ |x|=2$} & \mcol{1}{c}{$n=100,\ |x|=3$}\\
\mcol{1}{c|}{} & \mcol{1}{c|}{$\theta/\pi$} & \mcol{1}{c}{$\theta/\pi$}\\
[0.05cm]\hline
&&\\[-0.35cm]
$s_1,\, s_2$  & 0.12796 & 0.22172\\
$s_2,\, s_3$  & 0.05859 & 0.09844\\
$s_3,\, s_4$  & 0.03617 & 0.06070\\
$s_4,\, s_5$  & 0.02534 & 0.04264\\
$s_5,\, s_6$  & 0.01907 & 0.03220\\   
[.10cm]\hline
\end{tabular}
\end{center}
\end{table}

We illustrate this situation in Fig.~2 for the particular case $n=100$, $|x|=3$. Fig.~2(a) shows
the topology of the steepest descent paths $C_k$ through the saddles for $0.22172\pi<\theta\leq \fs\pi$. These paths with $k\leq 0$ are similar to those in Fig.~1, but the paths $C_k$ with $k\geq 1$ emanate from the logarithmic singularities $T_k$ and pass to infinity in the direction $\arg\,s\simeq\theta/(2\log\,|x|)$. The contribution to $\wp_n(x^{-2})$ in this case results only from the saddles with $k\leq 1$. Fig.~2(b) shows the critical value $\theta=0.22172\pi$ where the saddles $s_1$ and $s_2$ connect to produce a Stokes phenomenon. Fig.~2(c) illustrates the topology of the steepest descent paths for $\theta$ in the range $0.09844\pi<\theta<0.22172\pi$, where $\wp_n(x^{-2})$ picks up an additional contribution from the saddle $s_2$.
Finally, Fig.~2(d) shows the next critical value $\theta=0.09844\pi$ where the saddles $s_2$ and $s_3$ connect
to produce another Stokes phenomenon. 
If we denote by $K(\theta)$ the number of contributing saddles situated in Re$(s)>0$, then $K(\theta)=1$
for $0.22172\pi<\theta\leq\fs\pi$ and $K(\theta)=2$ for $0.09844\pi<\theta<0.22172\pi$.
This process of connecting pairs of saddles continues as $\theta$ approaches zero with the result that progressively more saddles contribute to the sum $\wp_n(x^{-2})$. When $\theta=0$, there are no such connections ($K(0)=\infty$) and all the saddles contribute to yield the result in (\ref{e28}).
\begin{figure}[ht]
	\begin{center}
		{\tiny($a$)}\includegraphics[width=0.35\textwidth]{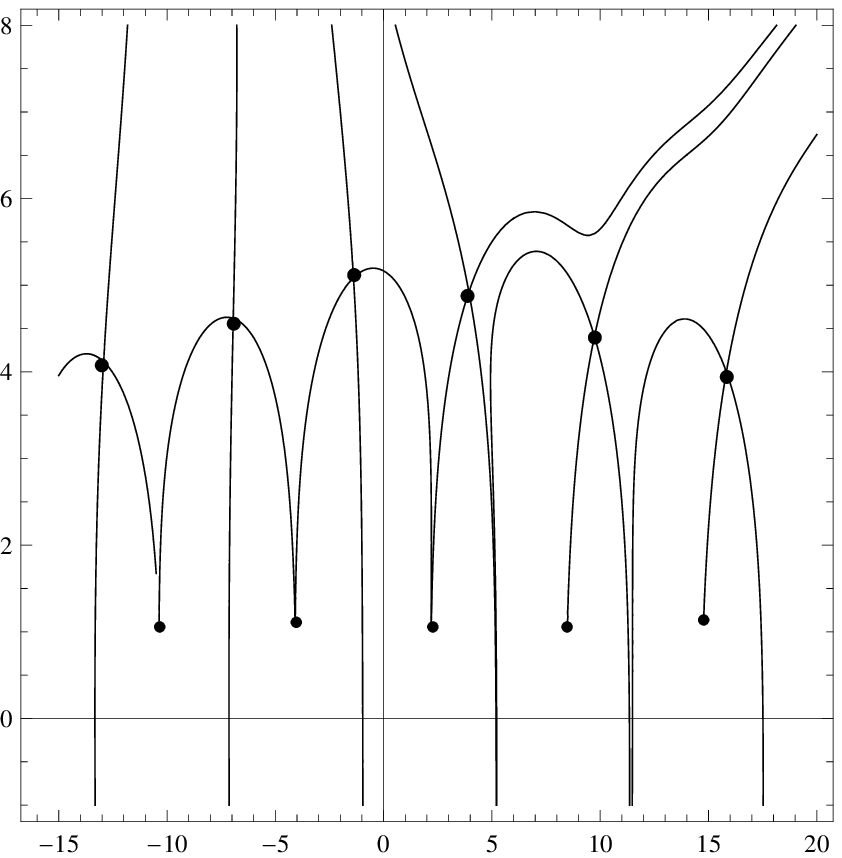}\hspace{1cm} {\tiny($b$)}\includegraphics[width=0.35\textwidth]{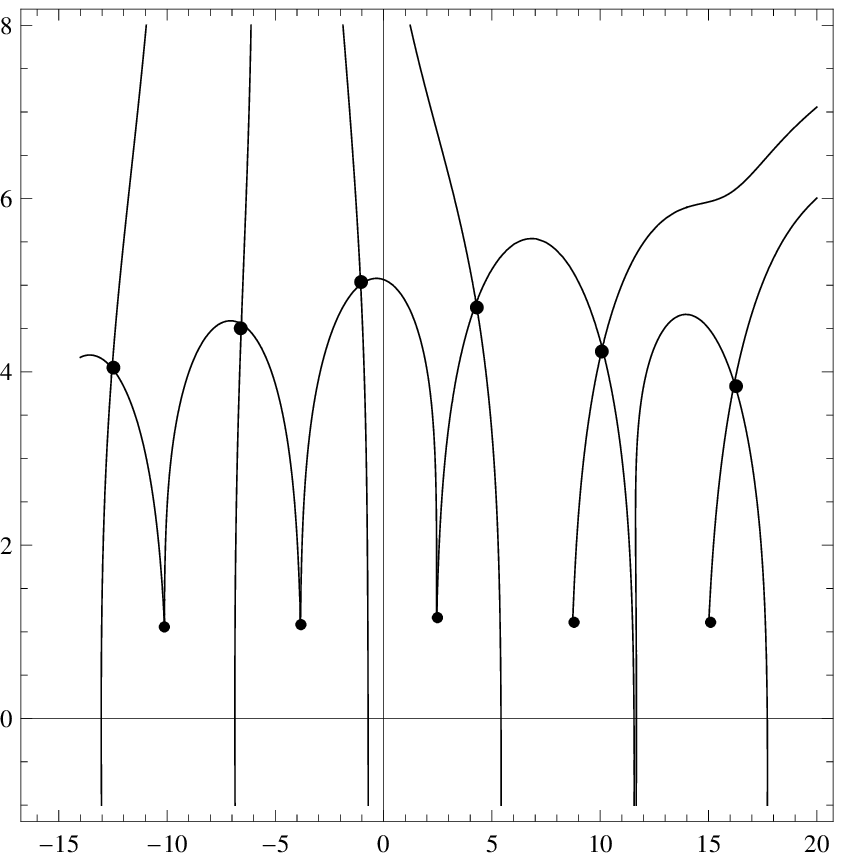}
		\vspace{0.6cm}
		
		{\tiny($c$)}\includegraphics[width=0.35\textwidth]{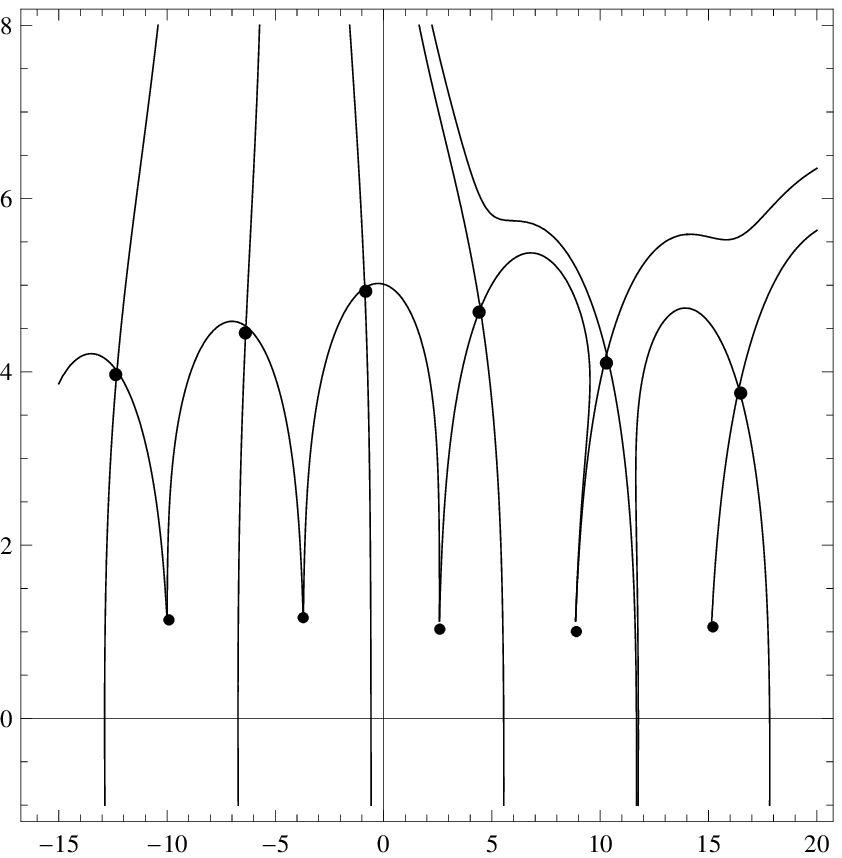}\hspace{1cm}
		{\tiny($d$)}\includegraphics[width=0.35\textwidth]{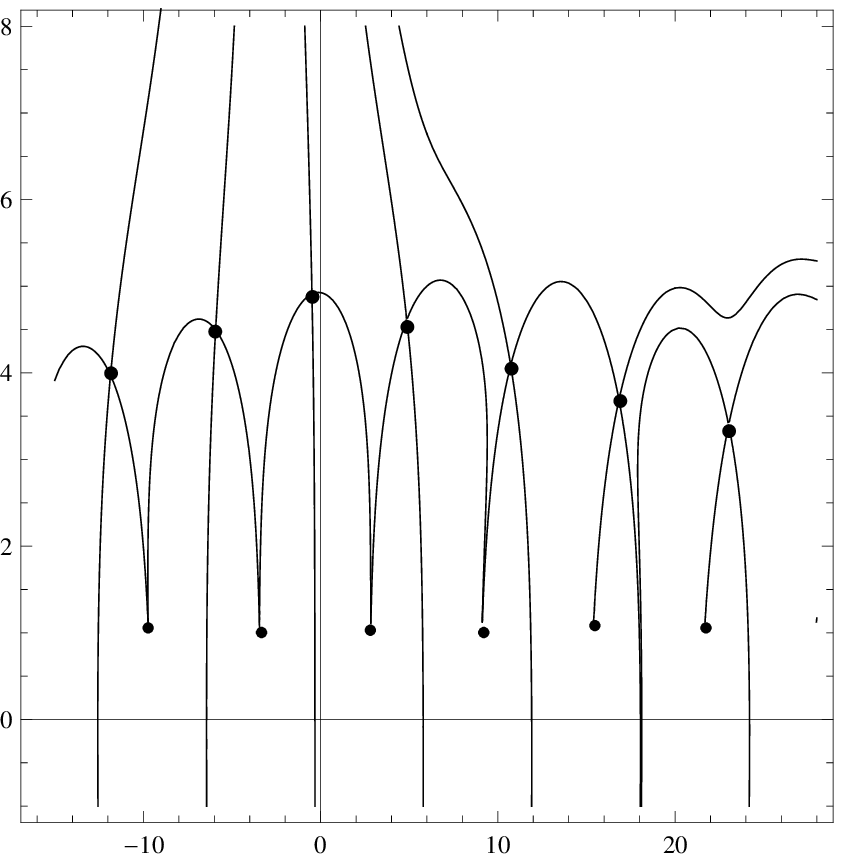}
	\caption{\small{The steepest descent paths ${\cal C}_k$ through the saddles $s_k$ (heavy dots) with $-2\leq k\leq 3$ (including $s_4$ in (d)) when $n=100$, $x=3e^{i\theta}$:
	(a) $\theta=0.30\pi$, (b) $\theta=0.22172\pi$, (c) $\theta=0.18\pi$, (d) $\theta=0.09844\pi$.
A Stokes phenomenon occurs in (b) and (d) where the two adjacent saddles $s_1, s_2$ and $s_2, s_3$ become connected. }}\label{f2}
	\end{center}
\end{figure}

In view of the above discussion, the expansion of $\wp_n(x^{-2})$ for complex $x$ when $|x|>1$ and $0\leq\theta\leq\fs\pi$ is given by
\bee\label{e41}
\wp_n(x^{-2})\sim \sum_{k=-\infty}^{K(\theta)}J_k ,
\ee
where the contributions $J_k$ have the expansions given in (\ref{e30a}) and (\ref{e30b}). The upper limit $K(\theta)$ varies with $\theta$: the precise transitions in $K(\theta)$ as $\theta$ decreases from $\fs\pi$ will depend on the values of $n$ and $|x|$.

\begin{table}[h]
\caption{\footnotesize{The absolute relative error in the computation of $\wp_n(x^{-2})$ for different $n$ and complex $x$ using (\ref{e41}).}}
\begin{center}
\begin{tabular}{l|c|c|c}
\mcol{1}{c|}{$\theta/\pi$} & \mcol{1}{c|}{$n=200,\ |x|=1.20$} & \mcol{1}{c|}{$n=200,\ |x|=1.50$} & \mcol{1}{c}{$n=1000,\ |x|=2.00$}\\
[0.05cm]\hline
&&&\\[-0.25cm]
0    & $5.919\times 10^{-9}$ & $5.329\times 10^{-7}$ & $2.206\times 10^{-6}$\\
0.10 & $3.391\times 10^{-7}$ & $1.303\times 10^{-6}$ & $3.564\times 10^{-6}$\\
0.20 & $4.315\times 10^{-7}$ & $3.413\times 10^{-6}$ & $3.216\times 10^{-6}$\\
0.30 & $3.856\times 10^{-7}$ & $1.077\times 10^{-5}$ & $6.450\times 10^{-6}$\\
0.40 & $2.808\times 10^{-7}$ & $7.425\times 10^{-6}$ & $1.029\times 10^{-5}$\\
0.50 & $5.379\times 10^{-8}$ & $4.089\times 10^{-6}$ & $2.245\times 10^{-6}$\\
[.10cm]\hline
\end{tabular}
\end{center}
\end{table}

\begin{figure}[t]
	\begin{center}	{\tiny($a$)}\includegraphics[width=0.35\textwidth]{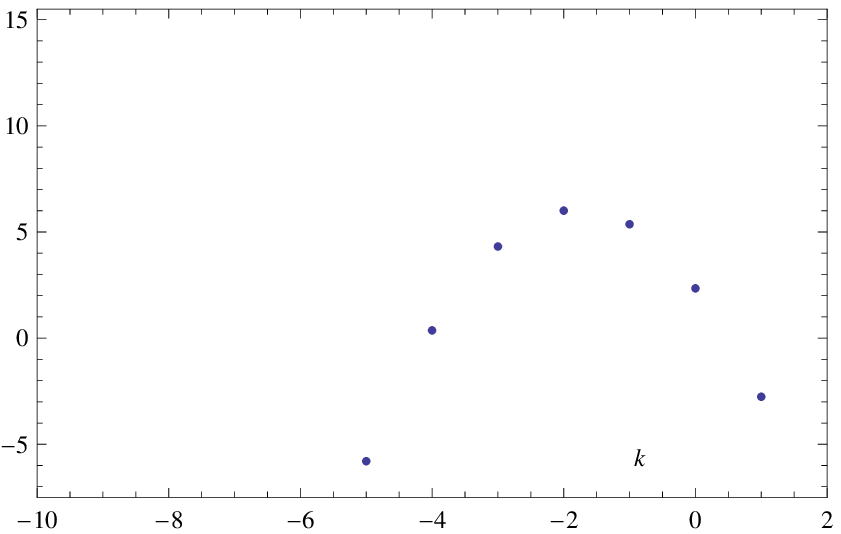}\qquad
	{\tiny($b$)}\includegraphics[width=0.35\textwidth]{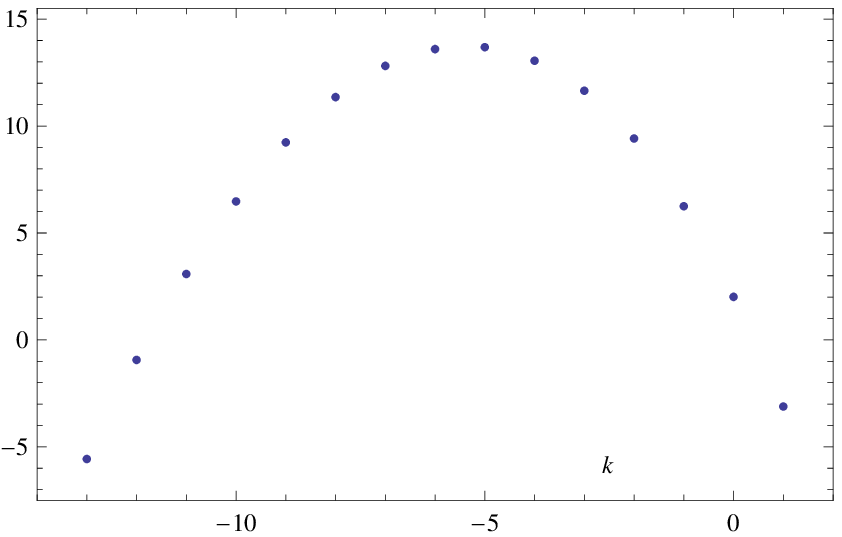}
\caption{\small{Values of $\log_{10}|J_k|$ for $k\leq 1$ when $n=200$ and $\theta=0.40\pi$: (a) $|x|=1.50$ and  (b) $|x|=1.10$.}}\label{f3}
	\end{center}
\end{figure}

In Table 5 we present the absolute relative error in the computation of $\wp_n(x^{-2})$ for complex $x$ for different $n$ and $|x|$ using (\ref{e41}). In each case we have employed the truncation index $j=3$ in the expansions of $J_k$. 
A feature of these calculations that is worthy of note is the location of the dominant saddles.
When $\theta=0$, the dominant saddle corresponds to $k=0$ and the saddles with $k=\pm 1$ yield a very small contribution. 
As $|x|\ra 1$ and $\theta\ra\fs\pi$ (that is, $|z|\ra 1$ and $\arg\,z\ra\pi$), the dominant saddles progressively shift to more negative $k$-values and, moreover, more saddles are found to make a significant contribution to $\wp_n(x^{-2})$. This last fact makes the calculation of $\wp_n(x^{-2})$ more difficult to determine in this limit.
We illustrate this feature in Fig.~3 where we display examples of the values of $\log_{10}|J_k|$ against the index $k$.

\vspace{0.6cm}

\begin{center}
{\bf Appendix: \ The quantities $Q_j(a_0,\lambda_0)$}
\end{center}
\setcounter{section}{1}
\setcounter{equation}{0}
\renewcommand{\theequation}{\Alph{section}.\arabic{equation}}
In this appendix we display the quantities $Q_j(a_0,\lambda_0)$ appearing in the coefficients $c_{j0}$ in (\ref{e2coeff}) for $1\leq j\leq 3$. These have the form
\begin{eqnarray}
Q_1(a_0,\lambda_0)&=&P_{11}(\lambda_0)-3a_0 P_{12}(\lambda_0),\nonumber\\
Q_2(a_0,\lambda_0)&=&P_{21}(\lambda_0)-6a_0 P_{22}(\lambda_0)+3a_0^2 P_{23}(\lambda_0)-48a_0^3 P_{24}(\lambda_0),\nonumber\\
Q_3(a_0,\lambda_0)&=&P_{31}(\lambda_0)+27a_0 P_{32}(\lambda_0)-9a_0^2 P_{33}(\lambda_0)+27a_0^3 P_{34}(\lambda_0)\nonumber\\
&&\hspace{4cm}-432a_0^4 P_{35}(\lambda_0)+4320a_0^5 P_{36}(\lambda_0),\label{a1}
\end{eqnarray}
where the $P_{jk}(\xi)$ ($1\leq k\leq 2j$) are polynomials in $\xi$ of degree $2j$ given by
\begin{eqnarray*}
P_{11}(\xi)&=&1+\xi+\xi^2,\qquad P_{12}(\xi)=1-4\xi+\xi^2,\\
\\
P_{21}(\xi)&=&(1+\xi+\xi^2)^2,\quad P_{22}(\xi)=13+5\xi-10\xi^2+5\xi^3+13\xi^4,\\
P_{23}(\xi)&=&67-328\xi+278\xi^2-328\xi^3+67\xi^4,\quad
P_{24}(\xi)=1-26\xi+66\xi^2-26\xi^3+\xi^4,\\
\\
P_{31}(\xi)&=&139+417\xi+402\xi^2+109\xi^3+402\xi^4+417\xi^5+139\xi^6,\\
P_{32}(\xi)&=&151+378\xi+308\xi^2+56\xi^3+308\xi^4+378\xi^5+151\xi^6,\\
P_{33}(\xi)&=&9271-3497\xi-10867\xi^2+766\xi^3-10867\xi^4-3497\xi^5+9271\xi^6,\\
P_{34}(\xi)&=&7349-48668\xi+45007\xi^2-24056\xi^3+45007\xi^4-48668\xi^5+7349\xi^6,\\
P_{35}(\xi)&=&203-5016\xi+18729\xi^2-24392\xi^3+18729\xi^4-5016\xi^5+203\xi^6,\\
P_{36}(\xi)&=&1-120\xi+1191\xi^2-2416\xi^3+1191\xi^4-120\xi^5+\xi^6.
\end{eqnarray*}

The coefficients $c_{jk}$ associated with the $k$th saddle follow from the expressions above and (\ref{e2coeff})
with $a_0$ and $\lambda_0$ replaced by $a_k$ and $\lambda_k$ defined in (\ref{e29}).
\vspace{0.6cm}

\noindent{\bf Acknowledgement}\ \ The author wishes to acknowledge access to the paper \cite{GN} before its publication and some helpful comments from the referees.


\vspace{0.6cm}

\begin{thebibliography}{99}
\footnotesize{
\bibitem{BDM}
J. BROWN, K. DILCHER AND D. MANNA, {\it Asymptotics of a sequence of sparse binomial-type polynomials}, Analysis {\bf 32} (2012) 231--245.

\bibitem{D}
R. B. DINGLE, {\it Asymptotic Expansions: Their Derivation and Interpretation}, Academic Press, London, 1973.

\bibitem{GN}
W. GAWRONSKI AND T. NEUSCHEL, {\it On a conjecture on sparse binomial-type polynomials by Brown, Dilcher and Manna}, Analysis and Applications, {\bf 12}(5) (2014) 511--522.

\bibitem{DLMF}
F. W. J. OLVER, D. W. LOZIER, R. F. BOISVERT AND C. W. CLARK (eds.), {\it NIST Handbook of Mathematical Functions}, Cambridge University Press, Cambridge, 2010.

\bibitem{P}
R. B. PARIS, {\it Hadamard Expansions and Hyperasymptotic Evaluation: An Extension of the Method of Steepest Descents}, Cambridge University Press, Cambridge, 2011.


}




\end{thebibliography}
\end{document}